\documentclass[%
 reprint,
 amsmath,amssymb,
 aps,
]{revtex4-2}
\usepackage{graphicx}
\usepackage{dcolumn}
\usepackage{bm}
\usepackage{array}
\usepackage[all,dvips]{xy}
\usepackage{epsfig}
\usepackage{soul}
\usepackage{color}
\usepackage{epstopdf}

\begin{document}

\preprint{APS/123-QED}

\title{Higher order interactions in Kuramoto oscillators with inertia}

\author{Patrycja Jaros$^{1,a}$, Subrata Ghosh$^{1,2}$, Dawid Dudkowski$^{1}$, Syamal K. Dana$^{1,3}$, Tomasz Kapitaniak$^1$}

\affiliation{$^{1}$Division of Dynamics, Lodz University of Technology, Stefanowskiego 1/15, 90-924 Lodz, Poland}

\affiliation{$^{2}$Center for Computational Natural Sciences and Bioinformatics, International Institute of Information Technology, Gachibowli, Hyderabad 500032, India.}

\affiliation{$^{3}$Department of Mathematics, National Institute of Technology,	Durgapur 713209, India}

\affiliation{$^{a}$ Author to whom correspondence should be addressed: patrycja.kuzma@p.lodz.pl}

\begin{abstract}
	
How higher-order interactions influence the dynamics of second order phase oscillators? We address this question using three coupled Kuramoto phase oscillators with inertia under both pairwise and higher order interactions, finding collective states, which are absent for pair-wise interactions. Chimera states appear in an expanded parameter region of the 3-node network by comparison to the network's response in the absence of the higher-order interactions. 
This small network of three nodes is investigated for different sets of phase lags and over a range of negative to positive coupling strength of pairwise as well as higher order or group interactions.  It demonstrates how important are the choices of the phase lag and the sign of the coupling for the emergent dynamics of the system.


 \end{abstract}
\maketitle


\section{INTRODUCTION}

 The collective behaviour of dynamical units in a network is an important topic
 of research in physics, biology, ecology, engineering, and social sciences \cite{Pikovsky_synchronization_book,Arenas_PhysRep2008, boccaletti2006complex,boccaletti2018synchronization}. Such studies help us gain understanding of many natural events such as the firing of fireflies, flocking of birds, swimming pool of fishes, behaviour of power grids and their failure, just to mention a few. The research in this direction started with the studies of coupled oscillators \cite{Pikovsky_synchronization_book} but in recent time it is more focused on collective behaviour of networks of dynamical units \cite{boccaletti2006complex,boccaletti2018synchronization,kuramoto2003chemical}. Synchronization is the main concern for understanding its stability and the causes of its destruction. Besides synchrony \cite{Pikovsky_synchronization_book}, chimera pattern \cite{abrams2004chimera,kuramoto2002coexistence,laing2009dynamics} is another important behaviour that evolves via self-organized symmetry breaking of synchronization; solitary state \cite{hellmann2020network} is also a matter of great concern that may lead to breakdown of synchrony such as failure of a power-grid. What has been started with a globally coupled network \cite{kaneko2015globally} and non-locally coupled ring of oscillators \cite{kuramoto2002coexistence}, it has now been extended to complex networks (small-world and scale free networks) \cite{boccaletti2006complex,boccaletti2018synchronization} with a variety of coupling forms (local coupling, synaptic coupling and others), and even to multilayer networks \cite{Zhang_PRL2015,anwar2022intralayer}. Initially all the efforts on network dynamics start with fixed network that do not change with time, the dynamics of the network are assumed to be restricted to the nodes. Later it has been extended to more realistic time-varying networks \cite{ghosh2022synchronized}, which assumes a change of the links of a network with time besides the dynamics of the nodes. However, the interaction between the nodes has been always assumed as pair-wise and restricted between any two nodes of a network at any instant of time. 
 \par A trend of research has started, in recent time, to go beyond the pairwise interactions and include non-pairwise group interactions, that is to include 3-nodes and 4-nodes higher order interactions and beyond \cite{battiston2021physics, ghorbanchian2021higher, chutani2021hysteresis}. Such higher order interactions between the nodes, which are projected to represent many real situations, reveal more information on the network dynamics, including existence of synchrony, chimeras and other collective behaviours that have been missing \cite{kovalenko2021contrarians,kundu2022higher,skardal2020higher,majhi2022dynamics,parastesh2022synchronization} in pair-wise interactive networks. 
It has been observed that, in several systems, functional brain dynamics, opinion formation \cite{Moussa_PlosOne2013}, disease spreading, random walk, ecological interaction, etc., cannot be adequately explained by pairwise interactions. Since these systems include many bodies or dynamical units, we must include higher order interactive terms in the dynamics to capture their behavior more comprehensively. A variety of new findings are thereby obtained in astrophysics \cite{stone2019statistical}, social science \cite{Vasilyeva_SciReport2021,de2020social}, ecology \cite{abrams1983arguments,guo2021evolutionary,chatterjee2022controlling} and other real world dynamical systems \cite{Patania_EPJ2017,Barbarossa_IEEE2020,Iacopini_NatCom2019,neuhauser2020multibody,wang2021simplicial,millan2021local}. Higher-order interaction can lead to the emergence of first order phase transitions \cite{skardal2020higher,matamalas2020abrupt} (asynchronous to synchronous state, onset of an outbreak of disease) in a system that was missing for pair-wise interaction of the same systems. Researchers have first started with applying higher-order interactions in network of phase oscillators. It is interesting to note that one study \cite{kovalenko2021contrarians} shows coupled phase oscillators exhibiting global synchronization under a combination of negative pairwise and positive higher-order interaction. Another study shows emergence of chimera \cite{kundu2022higher} that is possible  in a network of Kuramoto phase oscillators without inertia, which has been so far absent under the pairwise interactions. 
\par Motivated by such new behaviours, we plan to explore the effects of higher-order interaction in second order phase oscillators with inertia. Our study demonstrates that even a small number (three) of  second order phase oscillators can display several distinct states, when higher order interactions are taken into account and, which are so far missing for pair-wise interactions. We use three damped second order phase oscillators with diffusive pairwise interaction defined by the conventional sinusoidal coupling and in addition, higher order interactions of 2-simplices. The model and the results are presented in the following sections.

\section{THE MODEL}

We investigate Kuramoto model with inertia enriched by higher-order interactions for three coupled elements ($N=3$) that suffices to reveal rich dynamics of the system. The model is defined by the following equations:
  \begin{equation}
	\begin{aligned}
		m\ddot{\theta_{i}}+\varepsilon\dot{\theta_{i}}=\frac{\mu}{N} \sum_{j=1}^{N}sin(\theta_{j}-\theta_{i}-\alpha)\\
		+\frac{\gamma}{2N^2}\sum_{j=1}^{N}\sum_{k=1}^{N}\sin(\theta_{j}+\theta_{k}-2\theta_{i}-\alpha),    \label{eq1:} 
	\end{aligned}
\end{equation} 
where variable $\theta_i(t)$ describes the phase of the $i$-th oscillator, $i=1,2,3$. In our numerical simulations, we limit ourselves to a set of parameters, mass $m=1.0$ and damping  $\varepsilon=0.1$ when the phase oscillators remain excitable in isolation. The phase lag $\alpha$ and the coupling coefficients $\mu$ for pairwise interaction and $\gamma$ for higher order interactions are varied. This particular model (1) with  $\mu>0$ and $\gamma=0$ (no higher order interactions) has been widely studied earlier \cite{maistrenko2017smallest}.
We extend our investigations of the network dynamics with higher order interaction for both zero and non-zero lag parameter ($\alpha$).

 When there is no lag in the system ($\alpha=0$), we obtain stable fixed points for $\gamma > -3\mu$ and splay states for $\gamma > 3/2\mu$. For a range of negative $\mu$ and $\gamma$, there exist 2+1 antipodal fixed points, which we explain later while discussing the details in the next section. 
 Moving away from this parameter region, 2+1 phase locked states are obtained for larger $\mu$. The addition of a small lag ($\alpha=0.1$) induces a switching state between rotating waves and the chimera states. Synchrony still exist for $\gamma > -3\mu$. For a larger phase lag ($\alpha=1.6$), chimera states dominate, both for positive and negative $\gamma$ and $\mu$ values. As $\alpha$ crosses $\pi/2$ value,  we obtain complete synchrony for $\gamma < -3\mu$.
 Details of our numerical results are presented in the next sections.
  

 \section{Distinct states: ZERO PHASE LAG}
 
\begin{figure}[!ht]
  	\includegraphics[width=\linewidth]{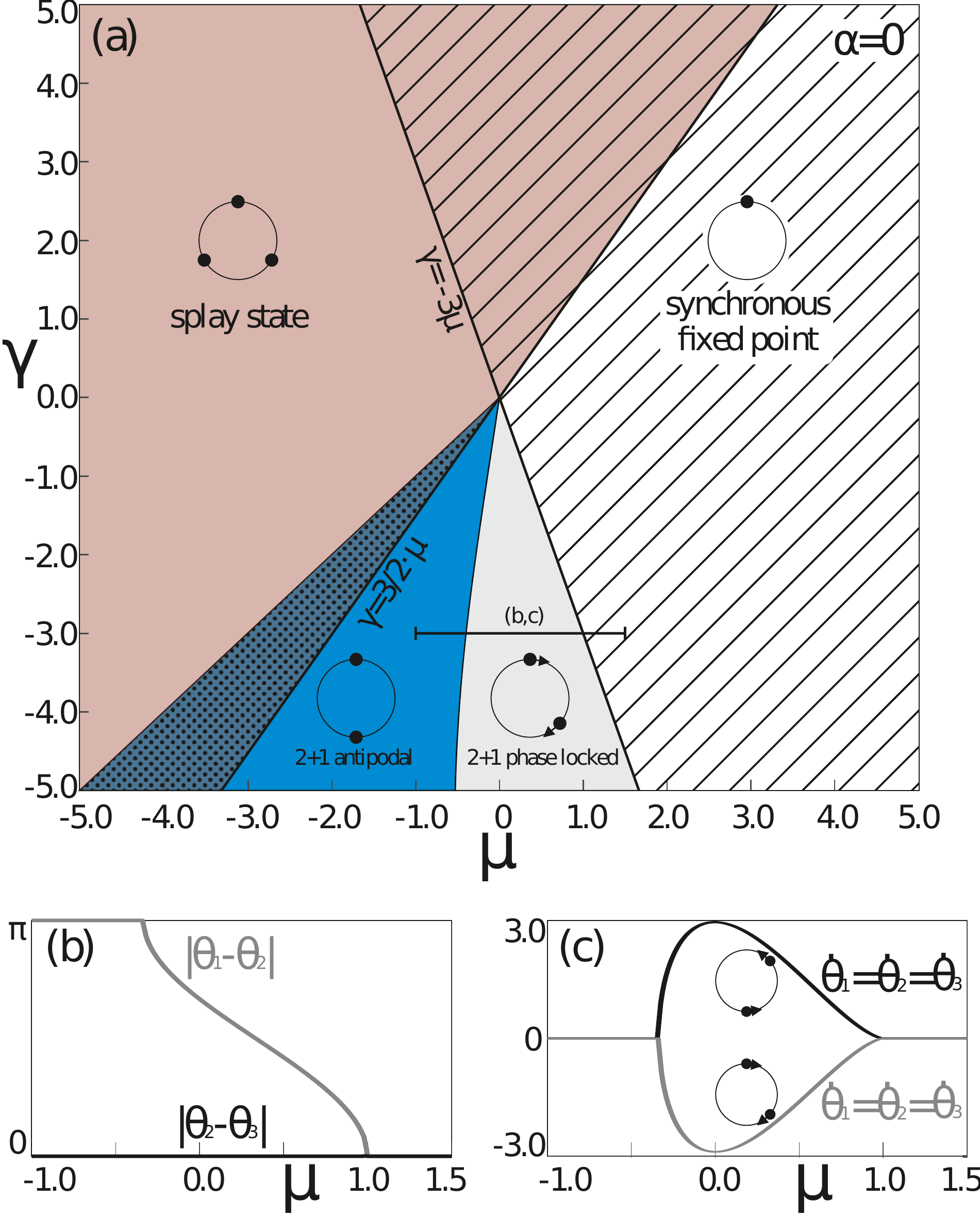}
  	\caption{(colour online) (a) $(\mu,\gamma)$-parameter plane for  $\alpha=0$, where hatched white area stands for synchronous fixed point,  hatched pink area for coexisting splay and synchronous state, pink for splay state, blue for 2+1 fixed antipodal point, light grey for 2+1 phase-locked state; (b, c) phase difference and phase velocity bifurcation diagram showing two coexisting 2+1 phase-locked states, respectively. Bifurcations are presented for $\mu$-interval indicated by a short horizontal line in (a). The circles with dots are schematic representations of each state.}
  	\label{fig1}
\end{figure}

 Firstly, we present the results of model (1) without a phase lag ($\alpha=0$). The results are summarized in a phase diagram of  ($\mu$, $\gamma$) parameter plane as shown in Fig.~1(a), where we see that without the higher order interactions ($\gamma=0)$, the network has two attractors: a synchronous fixed point (hatched region) for positive $\mu>0$ and splay state (pink region) for negative $\mu<0$; no motion is noticed.
 Non-zero $\gamma$ expands those regions in 2-dimensional space of ($\mu$, $\gamma$), where the synchronous fixed point is stable for $\gamma>-3\mu$ and a splay state exists for $\gamma>3/2\mu$, what have been confirmed analytically as well as numerically. These regions intersect only for certain positive values of $\gamma$ in the hatched pink region that indicates coexistence of those two solutions there.
 
 For negative $\gamma$, new states appear. In the blue region, 2+1 antipodal fixed point exists; it means that two elements are synchronized and one is not and it is located exactly at opposite phase to them (hence named as antipode). The antipodal fixed point, in particular, means all the elements have zero phase velocity. This state coexists with a splay state in a dotted blue region.  To the right of the blue region, in the light grey one, the only state with non-zero velocity can be found - the so called 2+1 phase-locked state. By this name, we mean that two elements are phase-synchronized, one is not, and they all move with the same constant velocity. As one may notice, the term antipodal is not used in this case, as the phase difference between the two clusters does not have to be equal to $\pi$.
 
 Although it may not be visible at first glance, the 2+1 phase-locked state is a very interesting behaviour. Higher-order interactions enable forcing of the system, which is why elements can rotate even with $\alpha=0$. Moreover, it appears that the 2+1 phase-locked state provides smooth transition from 2+1 fixed antipodal point to synchronous fixed point as shown in the bifurcation diagrams in Fig.~1(b, c)). In Fig.~1(b), the evolution of phase differences is shown. By setting $\gamma=-3$ and increasing $\mu$ from -1 to 1.5 (a short horizontal line marked in Fig.~1(a)), one may see how the phase difference between the synchronized elements $|\theta_2-\theta_3|$ stands at zero  and how the phase difference $|\theta_1-\theta_2|$ between a detached and one element from synchronized pair smoothly decreases, and finally also arriving at zero. When the  difference $|\theta_1-\theta_2|$ is equal to 0 or $\pi$, the system is in the equilibrium, but for any intermediate value $|\theta_1-\theta_2|\in(0,\pi)$,  the elements rotate, what can be seen in Fig.~1(c).  Fig.~1(c) also exhibits  that two symmetrical 2+1 phase-locked states coexist with each other (the black curve stands for counter-clockwise movement, while the grey one for clockwise rotation).
 
\section{Distinct states: SMALL PHASE LAG}
 Let us now examine how the dynamics of the network changes by the addition of a phase lag $\alpha$ in the coupling - firstly setting a small value of $\alpha=0.1$. The  system naturally gains a forcing that causes a movement: the splay state transforms into rotating wave, synchronous equilibrium starts synchronous rotations, 2+1 fixed antipodal point disappears - as elements start oscillatory motion, becoming 2+1 phase locked state.
 
 \begin{figure}
 \includegraphics[width=\linewidth]{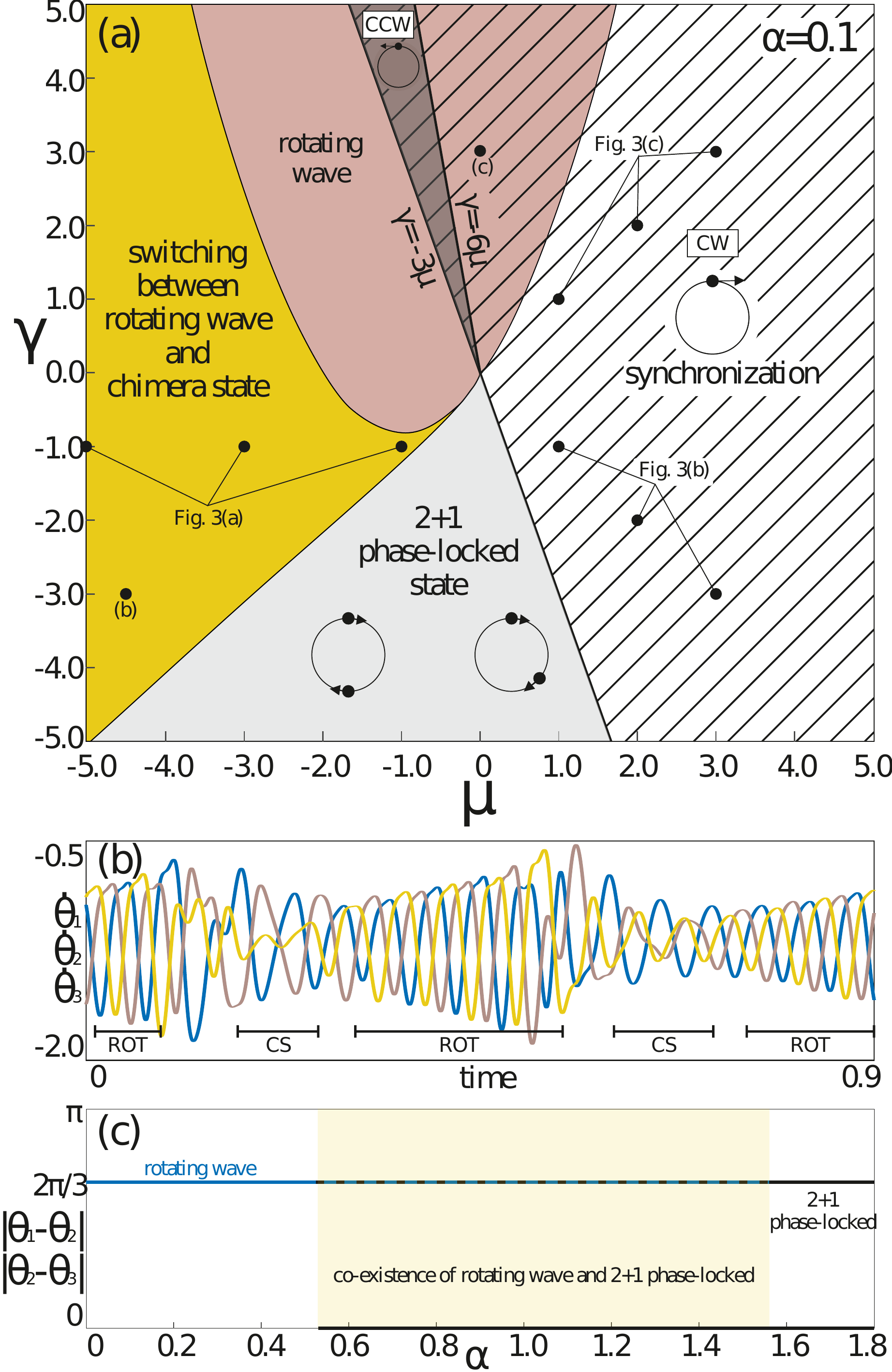}
 \caption{(colour online) (a) $(\mu,\gamma)$-parameter plane for  $\alpha=0.1$, where hatched area stands for synchronous rotations, pink for rotating wave, yellow for switching between rotating wave and chimera state, light grey for 2+1 phase-locked state; the shady region with label CW indicates clockwise and CCW counter-clockwise synchronous rotations; (b) velocity timeplot for $\mu=-4.5$, $\gamma=-3$ showing switching between rotating wave (ROT) and chimera state (CS); (c) bifurcation diagram for set $\mu=0$, $\gamma=3$ showing phase differences versus phase lag $\alpha$. The circles with dots are schematic representations of the states.}
 	\label{fig2}
 \end{figure} 
 
 It is worth mentioning that with careful selection of the coupling parameters and initial conditions, we can obtain clockwise or counter-clockwise synchronous rotations. As the equation for angular velocity of a synchronous state $\omega_s$ is given by: 
  \begin{center}
 $\omega_s=-\left( \dfrac{6\mu+\gamma}{9\varepsilon}\right)\sin\alpha $,
 \end{center} 
\noindent
 that means that counter-clockwise rotations exist for $\gamma>-6\mu$, what after taking into account a region of existence of synchronization, produces a grey, shady region indicated in Fig.~2(a). Choosing values of coupling  on the line $\gamma=-6\mu$, all oscillators are synchronized, but they do not move - their angular velocity  $\omega_s$ is equal to zero.
 What is worth mentioning here is the synchronous region that extends to a range of positive $\mu$ values (pairwise interaction), but spreads over a range of negative to positive $\gamma$ (groupwise interaction). This is similar to what has been reported earlier \cite{kovalenko2021contrarians}.
  In comparison to the case $\alpha=0$, the region of successor of the splay state rapidly shrinks, taking the shape similar to a parabola as shown in the phase diagram in the $(\mu,\gamma)$-parameter plane in Fig.~2(a). In the $(\mu,\gamma)$-parameter plane, a new player appears - {\it switching between rotating waves (ROT) and chimera states (CS)}, indicated by the yellow region. In this study, by chimera state we mean that two oscillators move with equal frequency, while the remaining one has different frequency. The switching between those two states (ROT and CS) means that the trajectory is once close to the rotating wave and after some time approaches the chimera state to return into the vicinity of the rotating wave. This scenario is being repeated infinitely. This kind of switching behavior has been already observed in the case of unidirectionally coupled Kuramoto model with inertia \cite{jaros2021chimera}. An exemplary velocity-time plot for the switching state is presented in Fig.~2(b) for $\mu=-4.5$, $\gamma=-3$, where we indicate the alternate time intervals of appearance of ROT and CS.
 
 Figure~2(c) presents the evolution of the system in a different manner: $(\mu,\gamma)$ is fixed and phase lag $\alpha$ is varied. For $\mu=0$, $\gamma=3$ and $\alpha=0$, splay state/rotating wave is obtained. When we make continuation by increasing $\alpha$ from $\alpha=0$ to  $\alpha=1.8$, we can see that the rotating wave loses it's stability for $\alpha\approx\pi/2$, and the system stabilizes at 2+1 phase-locked state with the phase difference equal to $2\pi/3$. However, if  $\alpha$ is decreased (starting from $\alpha$=1.8 to $\alpha=0$), it appears that the rotating wave coexists with 2+1 phase locked state for wide interval of $\alpha$ (the so-called hysteresis phenomenon).
 
 \begin{figure}[h]
	\includegraphics[width=\linewidth]{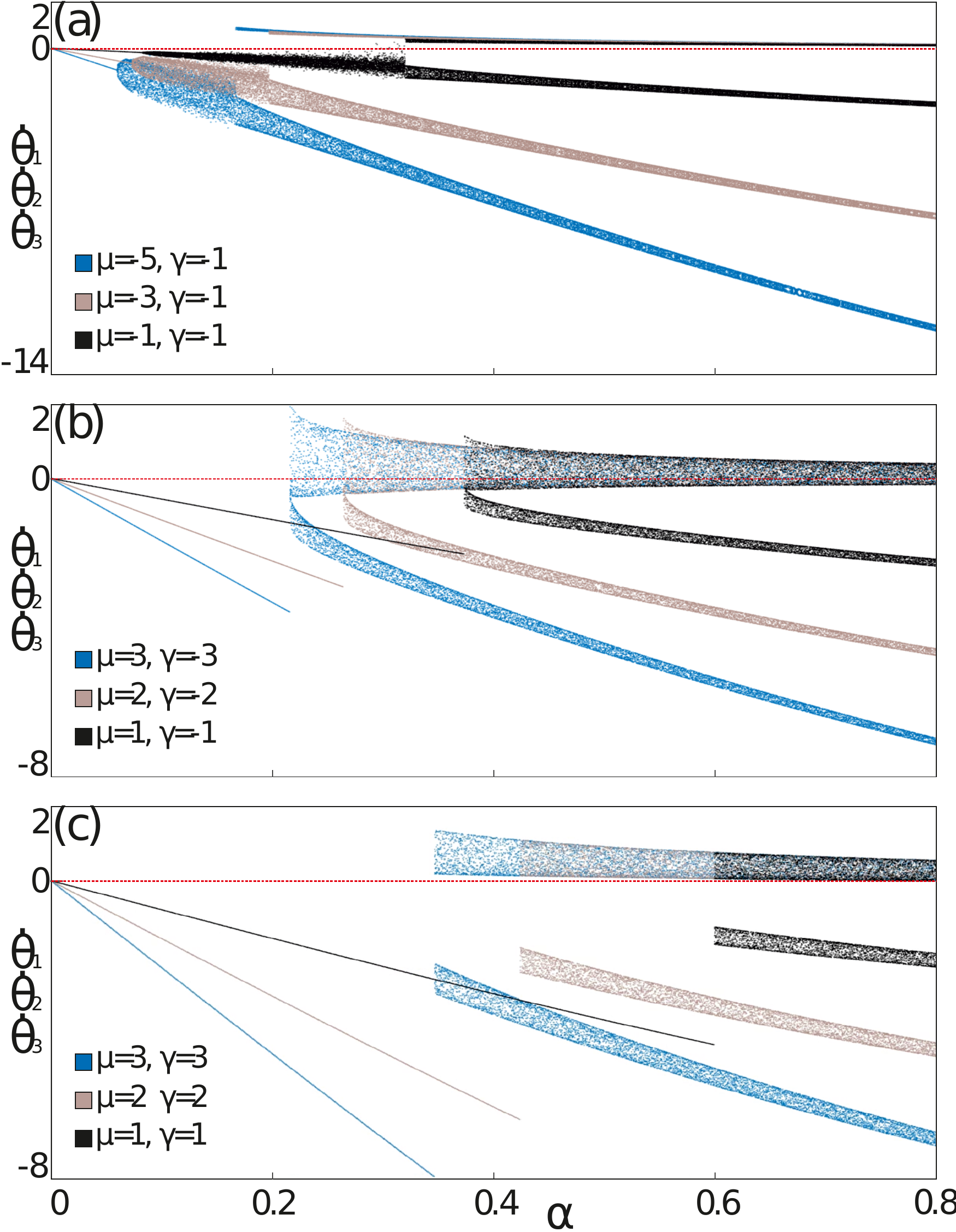}
	\caption{(colour online) Bifurcation of velocity against $\alpha$ for points chosen from Fig.~2(a). (a) Bifurcation diagram for increasing $\alpha$; (b,c) bifurcation for decreasing $\alpha$. The red line indicates zero velocity level.  } 
	\label{fig3}
\end{figure}

  \section{Distinct states: larger PHASE LAG }

 In Fig.~2(a), we indicate nine points (black circles with arms), for which bifurcation diagrams are presented in Fig.~3. The velocity-bifurcation diagrams show how chimera state (not only switching of it as for the case $\alpha=0.1$ in Fig.~2(a,c)) appears.
 
 Firstly, let us see what happens with the increase of $\alpha$, if we first pick up three points (black circles) of different parameter set of $(\mu,\gamma)$ in the yellow switching region (Fig.~2(a)), where $\mu<0$. Figure~3(a) presents how velocities of the oscillators evolve: in the beginning,  velocities of all the oscillators are negative and equal as the system is still in the rotating wave region.  When the latter one starts to shrink, a switching between the rotating wave and the chimera state takes place. In the bifurcation diagram, the switching is visible in the parabola-like shape with the vertex at the point of velocity of the rotating wave, which has been received as the last point before losing stability. With further increase of $\alpha$, one of the oscillators jumps to positive velocity level between 0.1 and 0.9 value  (depending on the set of chosen parameters), while the other two continues to evolve with a similar trend as for the rotating wave/chimera state switching (only the dispersion of velocity values is decreased). It is important to mention that if we create the bifurcation diagram here with decreasing $\alpha$, starting from $\alpha=0.8$, the diagram will be nearly identical (the coexistence of switching and chimera may occur in the very narrow $\alpha$-interval between those two states). Bifurcations are presented for three sets of parameters that enables us to conclude that the closer to the origin of all the regions (i.e. $(\mu,\gamma)$ = $(0,0)$), the greater $\alpha$ value needed for a chimera to appear.

 Now, let us move to a set of $(\mu,\gamma)$-pairs for positive $\mu>0$ and negative $\gamma<0$. This time the bifurcation diagrams are created for decreased $\alpha$. If we start from synchronous motion and increased $\alpha$, we would not obtain chimera, as synchronization is stable when $\alpha<\pi/2$ for those sets of chosen $(\mu,\gamma)$. In  Fig.~3(b) we see that there exists boundary after which the system jumps from chimera state to synchronous state (in this case common velocity does not mean rotating wave, but synchronization). There is no smooth transition as in case of Fig.~3(a). Again, we may see that if we take parameter point closer to $(\mu,\gamma)=(0,0)$, the longer (in the sense of increasing $\alpha$) chimera state does not occur. Totally analogously the situation looks  for positive $\mu>0$ and positive $\gamma>0$, which is presented in Fig.~3(c).  
  \begin{figure}
 \includegraphics[width=\linewidth]{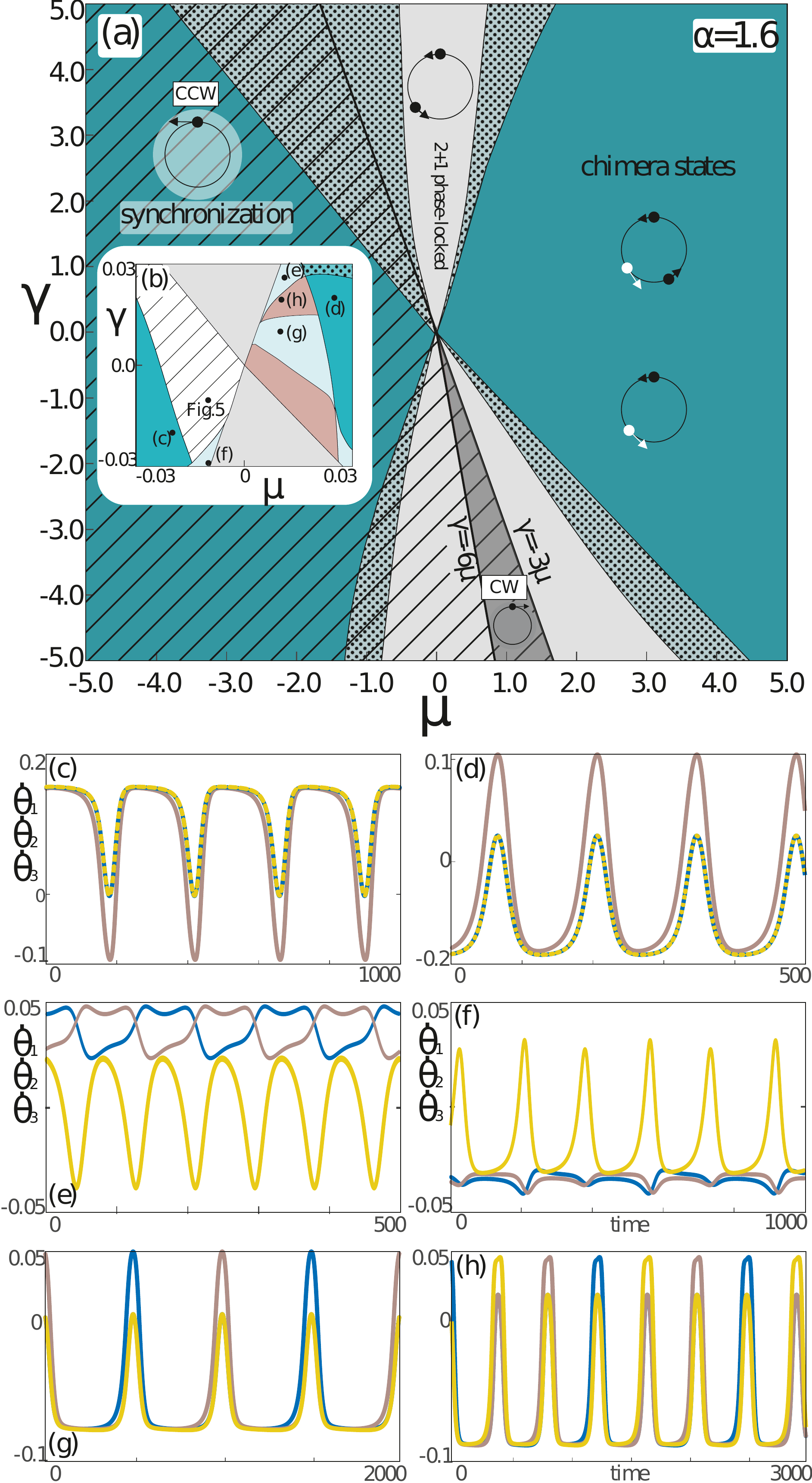}
 \caption{(colour online) (a) $(\mu,\gamma)$-parameter plane for  $\alpha=1.6$. Teal green colour stands for chimera states, light grey for 2+1 phase locked state, hatched region for synchronous rotations, dotted region for the coexistence of 2+1 phase locked and chimera state. The shady region with label CW indicates clockwise and CCW counter-clockwise synchronous rotations; the (b) inset is an enlarged version in the vicinity of $(\mu,\gamma)=(0,0)$: Regions in hatched white for globally stable synchronous rotation, dark green for in-phase chimera state, lighter green for anti-phase chimera state, pink colour for rotating waves. Exemplary time evolution plots of the variety of collective dynamics are presented in (c, d) in-phase chimera state, (e, f, g) anti-phase chimera state, and (h) rotating waves. The circles with dots are schematic representation of the movement for each state.}
  	\label{fig4}
  \end{figure} 
Figures~3(a-c) show that with the increase of $\alpha$, the size of the region of the occurrence of chimera states also increases. Now, we would like to see how the situation looks for $\alpha=1.6$, a value close to $\pi/2$, but a little larger. This case of $\alpha=1.6$ is equivalent to the case  $\alpha=\pi-1.6\approx1.5416$ as we consider both negative and positive values of the coupling. Results are qualitatively the same (only the direction of each movement is changed) if instead of $\mu$ and $\gamma$ we take $\tilde{\mu}=-\mu$ and $\tilde{\gamma}=-\gamma$. That is why there is no sense of further increase of $\alpha$ as we will obtain the same scenarios, only for horizontally and vertically flipped space.
Fig.~4(a) presents the $(\mu,\gamma)$-parameter plane for the case of $\alpha=1.6$.  Chimera state is already dominant for nearly all parameter plane. Other states that are present on the plane include synchronization (now stable for $\gamma<-3\mu$ as $\alpha>\pi/2$) and 2+1 phase-locked state. Often two of these states (or even three of them) coexist with each other (dotted region means coexistence of chimera state and phase-locked state, while additional hatching means coexistence of synchronous state). 

For $\alpha>\pi/2$, the region of synchronization is defined by $\gamma>-3\mu$. Generally, the performed synchronous rotations are counter-clockwise, apart from the grey, shady region defined by condition $\gamma>-6\mu$, where rotations change direction. As in the case of $\alpha=0.1$ for points from line  $\gamma=-6\mu$, the system is in the equilibrium.

As all the regions seem to originate from the $(\mu, \gamma)=(0,0)$ point, let us look at the enlarged version near this  point [see inset Fig.~4(b) in Fig.~4(a)]. In the  enlarged version, it can be seen that it is possible to set such parameters that synchronization is the only attractor in the system (white hatched region), what can be useful, if one searches for a safe region from the point of view of perturbations, which can occur in the system. This small globally stable complete synchronization regime do exists for purely repulsive interactions, both pairwise and groupwise.
\par Next, in the inset, we indicate different dynamics using two shades of green colour. Two different kinds of chimera states are shown: the in-phase (darker green) and the anti-phase (lighter green) one. The in-phase chimera state means that two frequency-synchronized elements are completely synchronized, while anti-phase chimera means that two elements with the same mean velocity are not in the phase (i.e. the phase difference do not have to be equal to $\pi$). There exist one more state, rotating waves, indicated by the pink regions.
\par Exemplary velocity-time plots of three oscillators in Fig.~4(c-h) show variety of possible states, which can be obtained just for small diversity of parameter values close to  $(\mu, \gamma)=(0,0)$.  Fig.~4(c,d) presents two similar in-phase chimera states, but showing opposite phase velocities. Anti-phase chimera states in Figs.~4(e-g) exhibit more diverse dynamics. In case (e) the desynchronized element's velocity is lower than the group performing anti-phase movement in velocity all the time.  In case (f)  the desynchronized element's velocity is higher than for the group with the same mean frequency, but not synchronized velocity trajectories. At first glance anti-phase chimera presented in Fig.~4(g) may be confused with rotating wave presented in Fig.~4(h). In the case of rotating wave from Fig.~4(h) we may see that the sequence of number of oscillator, which performs higher velocity jump takes the form: [1,2,3,1,2,3,...], while for chimera state Fig.~4(g) the same sequence consists of: [2,3,2,3,2,3...].

 \begin{figure}[h]
	\includegraphics[width=\linewidth]{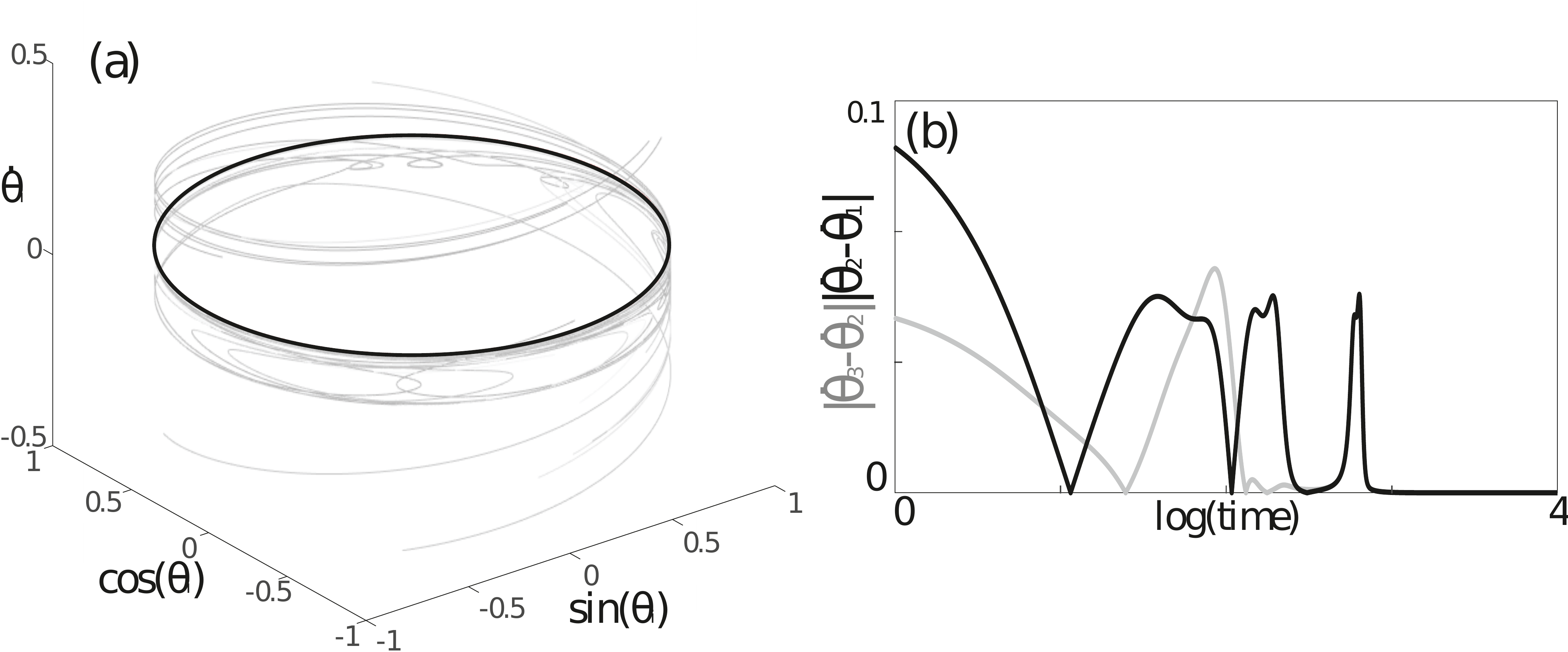}
	\caption{(colour online)  Trajectories for  $\alpha=1.6$, $\mu=-0.01$ and $\gamma=-0.01$ presented in (a) cylindrical space (grey curves stand for transient dynamics for few random initial conditions), black circle represents synchronized rotation of three units of the network and, (b) time plot of velocity differences $|\dot \theta_1-\dot \theta_2|$ and $|\dot \theta_2-\dot \theta_3|$; horizontal axis is in logarithmic scale when the time is varied from 1 to 10000.}
	\label{fig5}
\end{figure} 

The novel feature is the state of complete synchrony shown in Fig.~4(b) (see the inset, white hatched region), which is further exemplified here in Fig.~5 with the time evolution of three rotating elements of the network. It focusses on the parameter range,  $\alpha=1.6$, $\mu=-0.01$ and $\gamma=-0.01$, when both couplings are repulsive. According to Fig.~4(b), synchronization is the global attractor for such choice of parameters. Figure~5(a) shows  transient dynamics (grey lines) obtained for a few sets of random initial conditions. All of them converge to synchronized rotation, presented in a cylindrical space as a periodic orbit (black line). Fig.~5(b) shows plot of synchronization errors as $|\dot \theta_1-\dot \theta_2|$ and $|\dot \theta_2-\dot \theta_3|$ converging to zero in the long run.
	
 \section{DISCUSSIONS AND CONCLUSIONS}
  
 We have investigated the collective dynamics of a network of three nodes represented by the second order phase oscillators with inertia, popularly called Kuramoto oscillator with inertia. We set different values of phase lag and explored the parameter plane of the coupling terms to search for all possible collective states in the network dynamics.
  
 We have started our study with the simplest option - zero phase lag, when the additional coupling produces new collective behaviour: 2+1 antipodal point, where two elements are situated exactly opposite the third one and 2+1 phase-locked state - the only movement in the case of zero phase lag. The rotations of 2+1 phase-locked state is possible due to the negative $\gamma$, group interactions. Moreover, the 2+1 phase-locked state appears to show a smooth transition between 2+1 antipodal point and synchronization.  
  
 We have confirmed that the chimera states start to appear with the addition  of phase lag in the coupling. As the phase lag parameter is increased, the region of existence of chimera in parameter space expands. For $\alpha$ close to $\pi/2$, we have shown how big variety of different chimeras can be obtained for small variation of coupling parameters.
  
It seems that higher-order interactions support phase-locked states and coherence, in general. We have shown how the groupwise coupling affect ordinary synchronization. One may change the direction of motion and it's velocity by changing the coupling parameters, moreover, there exists a line in the $(\mu,\gamma)$ plane, where synchronization is restricted to the equilibrium state. Despite supporting coherence,  higher-order interactions do not kill chimera, but rather expands its territory in the parameter space, allowing it also to coexist with newly arisen phase-locked states.

\textbf{Acknowledgment}
 This work has been supported by the National Science Centre, Poland: OPUS Programme (Project No 2018/29/B/ST8/00457) and SONATA Programme (Project No 2019/35/D/ST8/00412).

\textbf{Data availability}: The data that support the findings of this study are available from the corresponding author upon request.


\begin{thebibliography}{70}
		
\bibitem{Pikovsky_synchronization_book} A. Pikovsky, M. Rosenblum, and J. Kurths, Synchronization: a universal concept in nonlinear sciences, Vol. 12 (Cambridge university press, 2003).

\bibitem{Arenas_PhysRep2008} A. Arenas, A. D´ıaz-Guilera, J. Kurths, Y. Moreno, and C. Zhou, Physics Reports 469, 93 (2008).

\bibitem{boccaletti2006complex} S. Boccaletti, V. Latora, Y. Moreno, M. Chavez, and D.-U. Hwang, Physics reports 424, 175 (2006).

\bibitem{boccaletti2018synchronization} S. Boccaletti, A. N. Pisarchik, C. I. Del Genio, and A. Amann, Synchronization: from coupled systems to complex networks (Cambridge University Press, 2018).

\bibitem{kuramoto2003chemical} Y. Kuramoto, “Chemical oscillations, waves and turbulence. mineola,” (2003).

\bibitem{abrams2004chimera} D. M. Abrams and S. H. Strogatz, Physical review letters 93, 174102 (2004).

\bibitem{kuramoto2002coexistence} Y. Kuramoto and D. Battogtokh, Nonlinear Phenomena in Complex Systems, 5, 380 (2002).

\bibitem{laing2009dynamics} C. R. Laing, Physica D: Nonlinear Phenomena 238, 1569 (2009).

\bibitem{hellmann2020network} F. Hellmann, P. Schultz, P. Jaros, R. Levchenko, T. Kapitaniak, J. Kurths, and Y. Maistrenko, Nature communications 11, 1 (2020).

\bibitem{kaneko2015globally}  K. Kaneko, Chaos: An Interdisciplinary Journal of Nonlinear Science 25, 097608 (2015).

\bibitem{Zhang_PRL2015} X. Zhang, S. Boccaletti, S. Guan, and Z. Liu, Physical Review Letters 114, 038701 (2015).

\bibitem{anwar2022intralayer} M. S. Anwar and D. Ghosh, Chaos: An Interdisciplinary Journal of Nonlinear Science 32, 033125 (2022).

\bibitem{ghosh2022synchronized} D. Ghosh, M. Frasca, A. Rizzo, S. Majhi, S. Rakshit, K. Alfaro-Bittner, and S. Boccaletti, Physics Reports 949, 1 (2022).

\bibitem{battiston2021physics} F. Battiston, E. Amico, A. Barrat, G. Bianconi, G. Ferraz de Arruda, B. Franceschiello, I. Iacopini, S. K{\'e}fi, V. Latora, Y. Moreno, et al., Nature Physics 17, 1093 (2021).

\bibitem{ghorbanchian2021higher} R. Ghorbanchian, J. G. Restrepo, J. J. Torres, and G. Bianconi, Communications Physics 4, 1 (2021).

\bibitem{chutani2021hysteresis} M. Chutani, B. Tadi{\'c}, and N. Gupte, Physical Review E 104, 034206 (2021).

\bibitem{kovalenko2021contrarians} K. Kovalenko, X. Dai, K. Alfaro-Bittner, A. Raigorodskii, M. Perc, and S. Boccaletti, Physical Review Letters 127, 258301 (2021).

\bibitem{kundu2022higher} S. Kundu and D. Ghosh, Physical Review E 105, L042202 (2022).

\bibitem{skardal2020higher} P. S. Skardal and A. Arenas, Communications Physics 3, 1 (2020).

\bibitem{majhi2022dynamics} S. Majhi, M. Perc, and D. Ghosh, Journal of the Royal Society Interface 19, 20220043 (2022).

\bibitem{parastesh2022synchronization} F. Parastesh, M. Mehrabbeik, K. Rajagopal, S. Jafari, and M. Perc, Chaos: An Interdisciplinary Journal of Nonlinear Science 32, 013125 (2022).

\bibitem{Moussa_PlosOne2013} M. Moussa{\"\i}d, J. E. K{\"a}mmer, P. P. Analytis, and H. Neth, PloS one 8, e78433 (2013).

\bibitem{stone2019statistical} N. C. Stone and N. W. Leigh, Nature 576, 406 (2019).

\bibitem{Vasilyeva_SciReport2021} E. Vasilyeva, A. Kozlov, K. Alfaro-Bittner, D. Musatov, A. Raigorodskii, M. Perc, and S. Boccaletti, Scientific Reports 11, 1 (2021).

\bibitem{de2020social} G. F. de Arruda, G. Petri, and Y. Moreno, Physical Review Research 2, 023032 (2020).

\bibitem{abrams1983arguments} P. A. Abrams, The American Naturalist 121, 887 (1983).

\bibitem{guo2021evolutionary} H. Guo, D. Jia, I. Sendi{\~n}a-Nadal, M. Zhang, Z. Wang, X. Li, K. Alfaro-Bittner, Y. Moreno, and S. Boccaletti, Chaos, Solitons \& Fractals 150, 111103 (2021).

\bibitem{chatterjee2022controlling} S. Chatterjee, S. Nag Chowdhury, D. Ghosh, and C. Hens, Chaos: An Interdisciplinary Journal of Nonlinear Science 32, 103122 (2022).

\bibitem{Patania_EPJ2017} A. Patania, G. Petri, and F. Vaccarino, EPJ Data Science 6, 18 (2017).

\bibitem{Barbarossa_IEEE2020} S. Barbarossa and S. Sardellitti, IEEE Transactions on Signal Processing 68, 2992 (2020).

\bibitem{Iacopini_NatCom2019} I. Iacopini, G. Petri, A. Barrat, and V. Latora, Nature Communications 10, 2485 (2019).

\bibitem{neuhauser2020multibody} L. Neuh{\"a}user, A. Mellor, and R. Lambiotte, Physical Review E 101, 032310 (2020).

\bibitem{wang2021simplicial} D. Wang, Y. Zhao, J. Luo, and H. Leng, Chaos: An Interdisciplinary Journal of Nonlinear Science 31, 053112 (2021).

\bibitem{millan2021local} A. P. Mill{\'a}n, R. Ghorbanchian, N. Defenu, F. Battiston, and G. Bianconi, Physical Review E 104, 054302 (2021).

\bibitem{matamalas2020abrupt} J. T. Matamalas, S. G´omez, and A. Arenas, Physical Review Research 2, 012049 (2020). 

\bibitem{maistrenko2017smallest} Y. Maistrenko, S. Brezetsky, P. Jaros, R. Levchenko, and T. Kapitaniak, Physical Review E 95, 010203 (2017).

\bibitem{jaros2021chimera} P. Jaros, R. Levchenko, T. Kapitaniak, and Y. Maistrenko, Chaos: An Interdisciplinary Journal of
Nonlinear Science 31, 103111 (2021).
	
\end{thebibliography}
\end{document}